\documentclass[10pt]{amsart}
\usepackage[cp1251]{inputenc}
\usepackage[english,russian]{babel}
\usepackage{amsmath}
\usepackage{longtable}
\usepackage{tikz}
\usepackage{amssymb}
\usepackage{amsfonts}

\newtheorem{proposition}{Proposition}

\newtheorem*{teoen}{Theorem}

\def\top
     {

     \begin{tabular}{c}
     \hphantom{aaaaaaaaaaaaaaaaaaaaaaaaaaaaaaaaaaaaaaaaaaaaaaaaaaaaaaaaaaaaaaaaaaaaaa} \\
     \end{tabular}

   {\flushleft
   \hspace{65mm}{\rm\small 
   } 
     \newline
         \hphantom{aaaaaaaaaaaaaaaaaaaaaaaaaaaaaaaaaaaaaaaaaaaaaaaaaaaaaaaa}{\rm\small MSC\ \ 20D60,~20D06 }
  }
}

\begin{document}
\renewcommand{\refname}{References}
\renewcommand{\proofname}{Proof.}
\thispagestyle{empty}

\title[Classification of Maximal Subgroups of Odd Index]{Classification of Maximal Subgroups of Odd Index \\ in Finite Simple Classical Groups: Addendum}
\author{{Natalia~V.~Maslova}}%
\address{Natalia~Vladimirovna~Maslova
\newline\hphantom{iii} Krasovskii Institute of Mathematics and Mechanics of the Ural Branch of the  Russian Academy of Science,
\newline\hphantom{iii} 16, S. Kovalevskaja St,
\newline\hphantom{iii} 620990, Ekaterinburg, Russia
\newline\hphantom{iii} Ural Federal University named after the first President of Russia B.~N.~ Yeltsin,
\newline\hphantom{iii} 19, Mira St,
\newline\hphantom{iii} 620002, Ekaterinburg, Russia}
\email{butterson@mail.ru}

\thanks{\rm The work is supported by the grant of the President of the Russian Federation for young scientists (grant no. MK-6118.2016.1) and by the Program for State Support of Leading Universities of the Russian Federation (agreement no. 02.A03.21.0006 of August 27, 2013).}

\top \vspace{1cm}
\maketitle {\small
\begin{quote}
\noindent{\sc Abstract. } A classification of maximal subgroups of odd index in finite simple groups was given by  Liebeck and Saxl and, independently, by Kantor in 1980s. In the cases of alternating groups or classical groups of Lie type over fields of odd characteristics, the classification was not complete.

The classification was completed by the author. In the cases of finite simple classical groups of Lie type we used results obtained in Kleidman's PhD thesis. However there is a number of  flaws in this PhD thesis. The flaws were corrected by Bray, Holt, and Roney--Dougal in 2013. In this note we provide a revision of the classification of maximal subgroups of odd index in finite simple classical groups over fields of odd characteristics.
\\

\noindent{\bf Keywords:} primitive permutation group, finite simple group, classical group, maximal subgroup, odd index.
 \end{quote}
}

\section{Introduction}

Liebeck and Saxl \cite{LiSa} and, independently, Kantor \cite{Kantor} have gaven a classification of finite primitive permutation groups of odd degree. It is one of the greatest results in the theory of finite permutation groups.

Both papers \cite{LiSa} and \cite{Kantor} contain lists of subgroups of finite simple groups that can be maximal subgroups of odd index. However, in the cases of alternating groups or classical groups of Lie type over fields of odd characteristics, neither in \cite{LiSa} nor in \cite{Kantor} it was described  which of the specified subgroups are precisely maximal subgroups of odd index. Thus, the problem of the complete classification of maximal subgroups of odd index in finite simple groups remained open.

The classification was completed by the author in \cite{Maslova1,Maslova2}. In \cite{Maslova1} we used results obtained by Kleidman \cite{Kleidman} and by Kleidman and Liebeck \cite{KlLi}. However there is a number of  flaws in Kleidman's PhD thesis \cite{Kleidman}. These flaws have been corrected in \cite{BrHoDo}. In this note we provide a revision of the classification of maximal subgroups of odd index in finite simple classical groups over fields of odd characteristics which was obtained in \cite{Maslova1}. In particular, we have made changes in statements of items $(6)$, $(10)$, and $(21)$ of \cite[Theorem~1]{Maslova1}.

\section{Terminology and Notation}

Throughout the paper we consider only finite groups, and thereby ''group'' means ''finite group''.

Our terminology and notation are mostly standard and can be found in \cite{BrHoDo,Atlas,KlLi,LiSa}.

\smallskip

If $G$ and $H$ are groups, $n$ is a positive integer, and $p$ is a prime, then we use notation $[n]$ for an arbitrary solvable group of order $n$, $\mathbb{Z}_n$ or just $n$ for the cyclic group of order $n$, $E_{p^n}$ or just $p^n$ for the elementary abelian $p$-group of order $p^n$, $G.H$ for an extension of $G$ by $H$, and $G:H$ for a split extension of $G$ by $H$.

Denote by $Soc(G)$ the {\it socle} of a group $G$ (i.\,e. the subgroup of $G$ generated by the non-trivial minimal normal subgroups of $G$). Recall that $G$ is {\it almost simple} if $Soc(G)$ is a non-abelian simple group.

\smallskip
The greatest integer power of $2$ dividing a positive integer $k$ is called the {\it $2$-part}
of $k$ and is denoted by $k_2$.

Let $m$ and $n$ be non-negative integers with $m=\sum_{i=0}^\infty a_i\cdot 2^i$ and $n=\sum_{i=0}^\infty b_i\cdot 2^i$, where $a_i,b_i\in\{0,1\}$. We write $m\preceq n$ if $a_i\leq b_i$ for every $i$.
It is clear that if $m\preceq n$, then $m \le n$. Moreover, $m\preceq n$ if and only if $n-m\preceq n$.

\medskip

\medskip

Let $q$ be a positive integer power of an odd prime $p$. Assume that $G$ is one of the simple classical
groups $PSL_n(q)$, $PSU_n(q)$, $PSp_n(q)$ for even $n$, $P\Omega_n(q)$ for odd $n$, or $P\Omega_n^\varepsilon(q)$ for even $n$, where
$\varepsilon \in \{+,-\}$. We will denote by $V$ the natural projective module of $G$, i.e. the vector space of dimension $n$ over a field $F$ with the corresponding bilinear form $\mathfrak{f}$ defined on this space, where
$F = \mathbb{F}_q$ for linear, symplectic, and orthogonal groups and $F = \mathbb{F}_{q^2}$ for unitary groups.
Note that if $\mathfrak{f}$ is non-degenerate, then for every non-degenerate subspace $U$ of $V$ we have (see \cite[Ch.~2]{KlLi})
\begin{equation}\label{DirectSum}V=U \oplus U^\bot.\end{equation}

\medskip

In the case of a non-degenerate symmetric bilinear form $\mathfrak{f}$ on $V$, the {\it discriminant} $D(V)$ of the corresponding quadratic form on $V$ is defined (see \cite[Ch.~2]{KlLi}). With a small abuse of terminology of \cite{KlLi}, we write $D(V)=1$ if $D(V)$ is a square in $F$, we write $D(V)=-1$ if $D(V)$ is a non-square in $F$.
In the case of the group $P\Omega_n^\varepsilon(q)$ for even $n$, the parameter $\varepsilon \in \{+,-\}$ is called
the {\it sign} of this group and of the corresponding vector space $V${\rm;} this parameter is denoted by $sgn(V)$.
For every non-degenerate subspace $U$ of even dimension $m$ from $V$,
the sign $\upsilon=sgn(U)$ of $U$ is defined, and the following formula is valid (see \cite[Proposition~2.5.10]{KlLi})

\begin{equation}\label{Discriminant}D(U)=D_m^\upsilon(q)=\begin{cases} 1 \mbox{ if } \upsilon=+ \mbox{ and }(q-1)m/4 \mbox{ is even}, \\  1 \mbox{ if } \upsilon=- \mbox{ and }(q-1)m/4 \mbox{ is odd}, \\ -1 \mbox{ otherwise.}\end{cases}\end{equation}
Moreover, it is known (see \cite[Proposition~2.5.11]{KlLi}) that if $U$ is a non-degenerate subspase of $V$, then \begin{equation}\label{MultRul1}D(V)=D(U)\cdot D(U^\bot), \end{equation} and  \begin{equation}\label{MultRul2}sqn(V)=sgn(U)\cdot sgn(U^\bot) \mbox{ if the dimensions of } U \mbox{ and } V \mbox{ are even}.\end{equation}

\medskip

Using the classification of finite simple groups, Aschbacher described in \cite{Asch} a large family of
natural geometrically defined subgroups of simple classical groups. He has subdivided this family
into eight classes $C_i$  for $1 \le i \le 8$, which are now called Aschbacher classes, and has proved that if a maximal subgroup of a simple classical group does not belong to the union of Aschbacher classes of the group, then this maximal subgroup is almost simple (for details see \cite{Asch}). The normal structure and the maximality of subgroups from Aschbacher classes of simple classical groups are known (see \cite{BrHoDo} for groups of dimension at most $12$ and \cite{KlLi} for groups of dimension at least $13$), and we will use these results.


\section{Main Result}

We prove the following theorem.

\begin{teoen}\label{Maslova_Revision} Let $G$ be one of the simple classical groups $PSL_n(q)$ for $n \ge 2$,
$PSU_n(q)$ for $n \ge 3$, $PSp_n(q)$ for even $n \ge 4$, $P\Omega_n(q)$ for odd $n \ge 7$, or $P\Omega^{\varepsilon}_n(q)$ for even $n \ge 8$ and $\varepsilon \in \{+,-\}${\rm;}
the number $q$ is always odd{\rm;} and let $V$ be the natural projective module of $G$. A subgroup $H$ of $G$ is maximal of odd index if and only if one of the following statements holds{\rm:}

$(1)$ $H = C_G(\sigma)$ for a field automorphism $\sigma$ of odd prime order $r$ of the group $G$, where $q=q_0^r${\rm;}

$(2)$ $G = PSL_n(q)$, $H$ is the stabilizer of a subspace of dimension $m$ of the space $V$, and $n \succeq  m${\rm;}

$(3)$ $G = PSU_n(q)$ or $G = PSp_n(q)$, $H$ is the stabilizer of a non-degenerate subspace of dimension $m$ of the space $V$, and  $n \succeq  m${\rm;}

$(4)$ $G = P\Omega_n(q)$, $H$ is the stabilizer of a non-degenerate subspace $U$ of even dimension $m$ of the space $V$, $D(U) = 1$, and $n \succeq  m${\rm;}

$(5)$ $G =P\Omega^\varepsilon_n(q)$, $H$ is the stabilizer of a non-degenerate subspace $U$ of dimension $m$ of the
space $V$, and one of the possibilities is realized:

$\mbox{ }\mbox{ }\mbox{ }\mbox{ }\mbox{ }$ $(i)$ $m$ is odd, $D(V)=-1$, and  $n-2 \succeq m-1${\rm;}

$\mbox{ }\mbox{ }\mbox{ }\mbox{ }\mbox{ }$ $(ii)$ $m$ is even, $D(U)=D(V )=-1$, $n-2 \succeq m-2$, $(q,m, sgn(U)) \not = (3, 2,+)$, and if $m=n/2$, then $\varepsilon=-${\rm;}

$\mbox{ }\mbox{ }\mbox{ }\mbox{ }\mbox{ }$ $(iii)$ $m$ is even, $D(U) = D(V ) = 1$, and  $n \succeq m${\rm;}

$(6)$ $G = PSL_n(q)$, $H$ is the stabilizer of a decomposition $V =\bigoplus V_i$ into a direct sum of
subspaces of equal dimension $m$, and one of the possibilities is realized{\rm:}

$\mbox{ }\mbox{ }\mbox{ }\mbox{ }\mbox{ }$ $(i)$ $m = 2^w \ge 2$ and $(n,m, q) \not= (4, 2, 3)${\rm;}

 $\mbox{ }\mbox{ }\mbox{ }\mbox{ }\mbox{ }$ $(ii)$ $m = 1$, $q \equiv 1
\pmod 4$, $(n, q) \not= (4, 5)$, and $q \ge 13$ in the case $n = 2${\rm;}

$(7)$ $G = PSU_n(q)$, $H$ is the stabilizer of an orthogonal decomposition $V =\bigoplus V_i$ into a direct sum of isometric subspaces $V_i$ of dimension $m$, and one of the possibilities is realized{\rm:}

 $\mbox{ }\mbox{ }\mbox{ }\mbox{ }\mbox{ }$ $(i)$ $m = 2^w \ge 2${\rm;}

 $\mbox{ }\mbox{ }\mbox{ }\mbox{ }\mbox{ }$ $(ii)$ $m = 1$, $q \equiv 3 \pmod 4$, and $(n, q) \not = (4, 3)${\rm;}

$(8)$ $G = PSp_n(q)$, $H$ is the stabilizer of an orthogonal decomposition $V =\bigoplus V_i$ into a direct
sum of isometric subspaces $V_i$ of dimension $m$, and $m = 2^w \ge 2${\rm;}

$(9)$ $G = P\Omega_n(q)$, $H$ is the stabilizer of an orthogonal decomposition $V = \bigoplus V_i$ into a direct
sum of isometric subspaces $V_i$ of dimension $1$, $q$ is prime, and $q \equiv \pm 3 \pmod 8${\rm;}

$(10)$ $G = P\Omega^\varepsilon_n(q)$, $H$ is the stabilizer of an orthogonal decomposition $V = \bigoplus V_i$ into a direct
sum of isometric subspaces $V_i$ of dimension $m$, $D(V)=1$, and one of the possibilities is realized{\rm:}

$\mbox{ }\mbox{ }\mbox{ }\mbox{ }\mbox{ }$ $(i)$ $m = 2^w \ge 2$, $D(V_i) = 1$, and $(m, q, sgn(V_i)) \not = (2, 3, -), (2, 5,+)${\rm;}

$\mbox{ }\mbox{ }\mbox{ }\mbox{ }\mbox{ }$ $(ii)$ $(n, \varepsilon) \not =(8,+)$, $m = 1$, $q$ is prime, and $q \equiv \pm 3 \pmod 8${\rm;}

$(11)$ $G = PSL_2(q)$ and $H \cong  PGL_2(q_0)$, where $q = q_0^2${\rm;}

$(12)$ $G = PSL_2(q)$ and $H \cong A_4$, where $q$ is prime and either $q=5$ or $q \equiv \pm 3, \pm 13 \pmod {40}${\rm;}

$(13)$ $G = PSL_2(q)$ and $H \cong S_4$, where $q$ is prime and $q \equiv \pm 7 \pmod {16}${\rm;}

$(14)$ $G = PSL_2(q)$ and $H \cong A_5$, where $q$ is prime and $q \equiv \pm 11, \pm 19 \pmod {40}${\rm;}

$(15)$ $G = PSL_2(q)$ and $H \cong D_{q+1}$, where $7 < q \equiv 3 \pmod 4${\rm;}

$(16)$ $G = PSU_3(5)$ and $H \cong M_{10}${\rm;}

$(17)$ $G = PSL_4(q)$ and $H \cong 2^4.A_6$, where $q$ is prime and $q \equiv 5 \pmod 8${\rm;}

$(18)$ $G = PSL_4(q)$ and $H \cong PSp_4(q).2$, where $q \equiv 3 \pmod 4${\rm;}

$(19)$ $G = PSU_4(q)$ and $H \cong 2^4.A_6$, where q is prime and $q \equiv 3 \pmod 8${\rm;}

$(20)$ $G = PSU_4(q)$ and $H \cong PSp_4(q).2$, where $q \equiv 1 \pmod 4${\rm;}

$(21)$ $G = PSp_4(q)$ and $H \cong 2^4.A_5$, where $q$ is prime and $q \equiv \pm 3 \pmod 8${\rm;}

$(22)$ $G = P\Omega_7(q)$ and $H \cong \Omega_7(2)$, where $q$ is prime and $q \equiv \pm 3 \pmod 8${\rm;}

$(23)$ $G = P\Omega^+_8(q)$ and $H \cong \Omega^+_8(2)$, where $q$ is prime and $q \equiv \pm 3 \pmod 8$.

\end{teoen}

\section{Proof of Theorem}

Using \cite{KlLi,Kantor} we conclude that the following proposition is valid.

\begin{proposition}\label{LSK_prop} Let $G$ be a simple classical group over a field of odd characteristic, $V$ be the natural projective module of $G$, and $H$ be a maximal subgroup of odd index of $G$.
Then one of the following statements holds{\rm:}

$(a)$ $H = N_G(C_G(\sigma))$ for a field automorphism $\sigma$ of  odd prime order $r$ of the group G{\rm;}

$(b)$ $H$ is the stabilizer of a non-degenerate {\rm(}arbitrary if $G=PSL_n(q)${\rm)} subspace of dimension $m$ of $V${\rm;}

$(c)$ $H$ is the stabilizer of an orthogonal {\rm(}arbitrary if $G=PSL_n(q)${\rm)} decomposition $V = \bigoplus\limits_{i=1}^{n/m} V_i$ with all $V_i$ isometric of dimension $m${\rm;}

$(d)$ $G$ is $P\Omega_7(q)$ or $P\Omega_8^+(q)$ and $H$ is $P\Omega_7(2)$ or $P\Omega_8^+(2)$, respectively, $q$ is prime,
and $q \equiv \pm 3 \pmod 8${\rm;}

$(e)$ $G = PSL_2(q)$ and $H$ is a dihedral group, $A_4$, $S_4$, $A_5$, or $PGL_2(q^{1/2})${\rm;}

$(f)$ $G = PSU_3(5)$ and $H \cong M_{10}$.

\end{proposition}

\smallskip

Note that results of \cite{LiSa} and \cite{Kantor} were formulated for orthogonal groups, but not for their isomorphic copies of other types. The maximal subgroups of groups $PSL_4(q)$, $PSU_4(q)$, and $PSp_4(q)$ are known \cite[Tables~8.8--8.13]{BrHoDo}, and following \cite{Maslova1} we consider these groups separately. However, it is not hard to obtain the similar results for the corresponding orthogonal groups.

\smallskip

Let $G$ be a simple classical group over a field of odd characteristic and $H$ be a subgroup of $G$ such that one of statements of Proposition \ref{LSK_prop} is valid.

\smallskip

The conditions when $G$ has a field automorphism of odd prime order are clear. If $H$ is a subgroup from the statement $(a)$ of Proposition \ref{LSK_prop}, then $C_G(\sigma)$ is a maximal subgroup of $G$ in view of \cite{BurGriLy}, and the index $|G:H|$ is odd in view of \cite{LiSa}. Thus, $H=C_G(\sigma)$, and the statement $(1)$ of Theorem holds.

Further consider simple linear, unitary, symplectic, and orthogonal groups separately.
We will use \cite[Tables~8.1-87]{BrHoDo} to check the maximality of subgroups in classical groups of dimension at most $12$ and \cite[Tables~3.5.A-F]{KlLi} to make the same for classical groups of dimension at least $13$.

\medskip

Assume that $G=PSL_n(q)$, where $n\ge 2$, $q$ is odd, and $(n,q)\not=(2,3)$.

If $H$ is the stabilizer of a subspace of dimension $m$ of $V$, then in view of \cite[Theorem~2]{Maslova1} the index $|G:H|$ is odd if and only if $n \succeq m$. In view of Tables~8.1, 8.3, 8.8, 8.18, 8.24, 8.35, 8.44, 8.54, 8.60, 8.70, and 8.76 of \cite{BrHoDo} and Table 3.5.A of \cite{KlLi} $H$ is maximal in $G$. Thus, the statement $(2)$ of Theorem holds.

\smallskip

If $H$ is the stabilizer of a decomposition $V = \bigoplus\limits_{i=1}^{n/m} V_i$ with all $V_i$ of dimension $m$,
then in view of \cite[Theorem~7]{Maslova1} the index $|G:H|$ is odd if and only if either $m=2^w\ge 2$, or $m=1$ and $q \equiv 1 \pmod 4$.
In view of Tables~8.1, 8.3, 8.8, 8.18, 8.24, 8.35, 8.44, 8.54, 8.60, 8.70, and 8.76 of \cite{BrHoDo} and Table 3.5.A of \cite{KlLi} $H$ is non-maximal in $G$ if and only if $n=2$ and $q \le 11$, or $m=1$ and $q=3$, or $(n,m,q)\in \{(4,1,5),(4,2,3)\}$. Thus, the statement $(6)$ of Theorem holds.

\smallskip

Consider additional abilities for the case $n=2$.

If $H \cong PGL_2(q_0)$ for $q=q_0^2$, then in view of  \cite[Tables~8.1,~8.7]{BrHoDo} $H$ is maximal in $G$, and it is easy to see that the index $|G:H|=(q_0^2+1)/2$ is odd. Thus, the statement $(11)$ of Theorem holds.

 Let $H \cong A_4$. In view of  \cite[Tables~8.1,~8.7]{BrHoDo} $G$ contains a maximal subgroup $H$ isomorphic to $A_4$ if and only if $q$ is prime and  either $q=5$ or $q \equiv  \pm 3, \pm 13 \pmod{ 40}$. It is clear that in this case the index $|G:H|$ is odd. Thus, the statement $(12)$ of Theorem holds.

 Let $H \cong S_4$. In view of \cite[Tables~8.1,~8.7]{BrHoDo} $G$ contains a maximal subgroup $H$ isomorphic to $S_4$ if and only if $q$ is prime and  $q \equiv  \pm 1 \pmod{8}$. It is clear that in this case the index $|G:H|$ is odd if and only if $q \equiv \pm 7 \pmod{16}$. Thus, the statement $(13)$ of Theorem holds.

 Let $H \cong A_5$. In view of  \cite[Tables~8.2,~8.7]{BrHoDo} $G$ contains a maximal subgroup $H$ isomorphic to $A_5$ if and only if either $q$ is prime and  $q \equiv \pm 1 \pmod{ 10}$, or $q=p^2$, where $p$ is a prime such that $p \equiv \pm 3 \pmod {10}$. It is clear that in the last case the index $|G:H|$ is even. In the first case the index $|G:H|$ is odd if and only if $q \equiv \pm 11, \pm 19 \pmod{40}$.
Thus, the statement $(14)$ of Theorem holds.

If $H$ is a dihedral group, then in view of  \cite[Tables~8.1,~8.7]{BrHoDo} either $H \cong D_{q-1}$, $H$ is the stabilizer of a decomposition $V = \bigoplus\limits_{i=1}^{2} V_i$ with both $V_i$ of dimension $1$, and the statement $(6)$ of Theorem holds; or $H \cong D_{q+1}$. In the last case $H$ is maximal in $G$ if and only if $q \not \in \{7,9\}$. It is clear that the index $|G:H|=q(q-1)/2$ is odd if and only if $q \equiv 3 \pmod 4$. Thus, the statement $(15)$ of Theorem holds.

\smallskip

Consider the case $n=4$. In view of the Aschbacher theorem \cite{Asch} and \cite[Tables~8.8,~8.9]{BrHoDo} any maximal subgroup of $G$ is either contained in $\bigcup_{i \in \{1,2,3,5,6,8\}}C_i(G)$ or is almost simple. Let $H$ be a subgroup of $G$.

Assume that $H \in C_1(G)$ and $H$ is maximal in $G$. Then in view of \cite[Table~8.8]{BrHoDo} $H$ is the stabilizer of a subspace of $V$, and the statement $(2)$ of Theorem holds.

Assume that  $H \in C_2(G)$ and $H$ is maximal in $G$. Then in view of \cite[Table~8.8]{BrHoDo} $H$ is the stabilizer of a decomposition $V = \bigoplus\limits_{i=1}^{n/m} V_i$ with all $V_i$ of dimension $m$, and the statement $(6)$ of Theorem holds.

Assume that $H \in C_3(G)$. Then in  view of \cite[Tables~8.8]{BrHoDo} the preimage of $H$ in $SL_4(q)$ is isomorphic to $SL_2(q^2).\mathbb{Z}_{q+1}.\mathbb{Z}_2$ and $$|G:H|_2=\frac{(q^4-1)_2(q^3-1)_2(q^2-1)_2}{2(q^4-1)_2(q+1)_2}\ge 2.$$

Assume that $H \in C_5(G)$. Then $H = C_G(\sigma)$ for a field automorphism $\sigma$ of prime order $r$ of the group $G$.
If $r$ is odd, then the statement $(1)$ of Theorem holds. If $r=2$, then in view of \cite[Table~8.8]{BrHoDo} the preimage of $H$ in $SL_4(q)$ is isomorphic to $SL_4(q_0).[(q_0+1,4)]$ for $q=q_0^2$ and
$$|G:H|_2=\frac{(q^4-1)_2(q^3-1_2)(q^2-1)_2}{(q_0^4-1)_2(q_0^3-1_2)(q_0^2-1)_2(q_0+1,4)}\ge 2.$$

In view of \cite[Table~8.8]{BrHoDo} $G$ contains a maximal subgroup $H \in C_6(G)$ if and only if $q$ is prime and $q \equiv 1 \pmod 4$. In view of \cite[Proposition~4.6.6]{KlLi}

$$H \cong \begin{cases} 2^4.S_6 \mbox{ if } q \equiv 1 \pmod 8,\\ 2^4.A_6 \mbox{ if } q \equiv 5 \pmod 8. \end{cases}$$

Note that we have $$|G:H|_2=\begin{cases}\frac{(q^4-1)_2(q^3-1)_2(q^2-1)_2}{2^{10} } > 2 \mbox{ if } q \equiv 1 \pmod 8,\\
\frac{(q^4-1)_2(q^3-1)_2(q^2-1)_2}{2^{9} }=1 \mbox{ if } q \equiv 5 \pmod 8. \end{cases}$$

Thus, the statement $(17)$ of Theorem holds.

Assume that $H \in C_8(G)$. Then view of \cite[Table~8.8]{BrHoDo} and \cite[Proposition~4.8.3]{KlLi} either $H$ is isomorphic to $PSp_4(q).[\frac{(q-1,2)^2}{(q-1,4)}]$, or the preimage of $H$ in $SL_4(q)$ is isomorphic to one of the following groups: $SO_4^\varepsilon(q).[(q-1,4)]$ for $\varepsilon \in\{+,-\}$ or $SU_4(q_0).\mathbb{Z}_{(q_0-1,4)}$, where $q=q_0^2$.

Then view of \cite[Table~8.8]{BrHoDo} $G$ always contains a maximal subgroup $H \cong PSp_4(q).\left[\frac{(q-1,2)^2}{(q-1,4)}\right]$ and $$|G:H|_2=\frac{2(q^4-1)_2(q^3-1)_2(q^2-1)_2(q-1,4)}{(q^4-1)_2(q^2-1)_2(q-1,4)(q-1,2)^2}=\frac{(q-1)_2}{2}.$$

Now it is easy to see that $|G:H|_2=1$ if and only if $q \equiv 3 \pmod 4$. Thus, $H \cong PSp_4(q).\mathbb{Z}_2$, and the statement $(18)$ of Theorem holds.

If the preimage of $H$ in $SL_4(q)$ is isomorphic to $SO_4^\varepsilon(q).[(q-1,4)]$ for  $\varepsilon \in\{+,-\}$, then $$|G:H|_2=\frac{(q^4-1)_2(q^3-1)_2(q^2-1)_2}{(q-1,4)(q^2-\varepsilon1)_2(q^2-1)_2}\ge \frac{(q^4-1)_2}{(q^2-\varepsilon1)_2}\ge 2.$$

If $q=q_0^2$ and the preimage of $H$ in $SL_4(q)$ is isomorphic to  $SU_4(q_0).\mathbb{Z}_{(q_0-1,4)}$, then $$|G:H|_2=\frac{(q^4-1)_2(q^3-1)_2(q^2-1)_2}{(q_0^4-1)_2(q_0^3+1)_2(q_0^2-1)_2(q_0-1,4)}\ge
\frac{(q_0^4+1)_2(q_0^3-1)_2(q_0^2+1)_2}{4}\ge 2$$

Finally, assume that $H$ is a maximal almost simple subgroup of $G$. Then in view of \cite[Table~8.9]{BrHoDo} $H \in \{PSL_2(7), A_7, PSU_4(2)\}$. So, $|H|_2 \le 2^6$ in view of \cite{Atlas}, and it is easy to see that $$|G|_2=\frac{(q^4-1)_2(q^3-1)_2(q^2-1)_2}{(q-1,4)}\ge 2^7.$$.

Thus, the index $|G:H|$ is even.

\medskip

Assume that $G=PSU_n(q)$, where $n\ge 3$ and $q$ is odd.

 Let $H$ be the stabilizer of a non-degenerate subspace $U$ of dimension $m$ of $V$. Note that in view of (\ref{DirectSum}) $H$ is the stabilizer of $U^\bot$, and we can assume that $m < n/2$ (if $m=n/2$, then $H$ is contained in a subgroup which will be considered in the next paragraph). In view of \cite[Theorem~3]{Maslova1} the index $|G:H|$ is odd if and only if $n \succeq m$. Moreover, if $n \succeq m$ then $m \not = n/2$. If $m \not = n/2$, then in view of Tables~8.5, 8.10, 8.20, 8.26, 8.37, 8.46, 8.56, 8.62, 8.72, and 8.78 of \cite{BrHoDo} and Table 3.5.B of \cite{KlLi} $H$ is maximal in $G$. Thus, the statement $(3)$ of Theorem holds for unitary groups.

If $H$ is the stabilizer of an orthogonal decomposition $V = \bigoplus\limits_{i=1}^{n/m} V_i$ with all $V_i$ isometric of dimension $m$, then in view of \cite[Theorem~8]{Maslova1} the index $|G:H|$ is odd if and only if either $m=2^w\ge 2$, or $m=1$ and $q \equiv 3 \pmod 4$. In view of Tables~8.5, 8.10, 8.20, 8.26, 8.37, 8.46, 8.56, 8.62, 8.72, and 8.78 of \cite{BrHoDo} and Table 3.5.B of \cite{KlLi} $H$ is non-maximal if and only if $m=1$  and $(n,q) \in \{(3,5), (4,3)\}$. Thus, the statement $(7)$ of Theorem holds.

If $G = PSU_3(5)$ and $H \cong M_{10}$, then in view of \cite{Atlas} $H$ is a maximal subgroup of odd index in $G$. Thus, the statement $(16)$ of Theorem holds.

Consider the case $n=4$. In view of the Aschbacher theorem \cite{Asch} and \cite[Tables~8.10,~8.11]{BrHoDo} any maximal subgroup of $G$ is either contained in $\bigcup_{i \in \{1,2,5,6\}}C_i(G)$ or is almost simple. Let $H$ be a subgroup of $G$.

Assume that $H \in C_1(G)$. Then  in view of \cite[Table~8.10]{BrHoDo} either $H$ is the stabilizer of a non-degenerate subspace of $V$, and the statement $(3)$ of Theorem holds, or $H$ is the stabilizer of a totally singular subspace of $V$ of dimension $m$. In the last case $1 \le m \le 2$, and in view of \cite[Table~8.10]{BrHoDo} the preimage of $H$ in $SU_4(q)$ is isomorphic to $E_q^{1+4}.SU_2(q):\mathbb{Z}_{q^2-1}$ for $m=1$ or to $E_q^4:SL_2(q^2):\mathbb{Z}_{q-1}$ for $m=2$, and $$|G:H|_2=\begin{cases}\frac{(q^4-1)_2(q^3+1)_2(q^2-1)_2}{(q^2-1)_2^2 }\ge 2 \mbox{ for }m=1,\\\frac{(q^4-1)_2(q^3+1)_2(q^2-1)_2}{(q^4-1)_2(q-1)_2} \ge 2 \mbox{ for }m=2.\end{cases}$$

Assume that  $H \in C_2(G)$. Then in view of \cite[Table~8.10]{BrHoDo} either $H$ is the stabilizer of an orthogonal decomposition $V = \bigoplus\limits_{i=1}^{n/m} V_i$ with all $V_i$ isometric of dimension $m$, and the statement $(7)$ of Theorem holds, or $H$ is the stabilizer of a decomposition $V=V_1\oplus V_2$, where $V_1$ and $V_2$ are totally singular subspaces of $V$ of dimension $2$. In the last case in view of \cite[Table~8.10]{BrHoDo} the preimage of $H$ in $SU_4(q)$ is isomorphic to $SL_2(q^2).\mathbb{Z}_{q-1}.\mathbb{Z}_2$ and  $$|G:H|_2=\frac{(q^4-1)_2(q^3+1)_2(q^2-1)_2}{2(q^4-1)_2(q-1)_2} \ge 2.$$

Assume that $H \in C_5(G)$. Then in view of \cite[Table~8.10]{BrHoDo} either $H = C_G(\sigma)$ for a field automorphism $\sigma$ of prime odd order of the group $G$, and the statement $(1)$ of Theorem holds, or in view of \cite[Table~8.10]{BrHoDo} and \cite[Proposition~4.5.6]{KlLi} either $H$ is isomorphic to $PSp_4(q).\left[\frac{(q+1,2)(q-1,2)}{(q+1,4)}\right]$, or the preimage of $H$ in $SU_4(q)$ is isomorphic to $SO_4^\varepsilon(q).[(q+1,4)]$ for $\varepsilon \in\{+,-\}$.

In view of \cite[Table~8.10]{BrHoDo} $G$ always contains a maximal subgroup $$H \cong PSp_4(q).\left[\frac{(q+1,2)(q-1,2)}{(q+1,4)}\right],$$ and
$$|G:H|_2=\frac{2(q^4-1)_2(q^3+1)_2(q^2-1)_2(q+1,4)}{(q+1,4)(q^4-1)_2(q^2-1)_2(q+1,2)(q-1,2)}=\frac{(q+1)_2}{2}.$$

Note that $|G:H|$ is odd if and only if $q \equiv 1 \pmod 4$. Thus, $H \cong PSp_4(q).\mathbb{Z}_2$, and the statement $(20)$ of Theorem holds.

If the preimage of $H$ in $SU_4(q)$ is isomorphic to $SO_4^\varepsilon(q).[(q+1,4)]$ for $\varepsilon \in\{+,-\}$, then
$$|G:H|_2=\frac{(q^4-1)_2(q^3+1)_2(q^2-1)_2}{(q+1,4)(q^2-\varepsilon1)_2(q^2-1)_2}\ge \frac{(q^4-1)_2}{(q^2-\varepsilon1)_2}\ge 2.$$

In view of \cite[Table~8.10]{BrHoDo} $G$ contains a maximal subgroup $H\in C_6(G)$ if and only if $q$ is prime and $q \equiv 3 \pmod 4$. In view of \cite[Proposition~4.6.6]{KlLi}

$$H \cong \begin{cases} 2^4.S_6 \mbox{ if } q \equiv 7 \pmod 8,\\ 2^4.A_6 \mbox{ if } q \equiv 3 \pmod 8. \end{cases}$$

Note that we have $$|G:H|_2=\begin{cases}\frac{(q^4-1)_2(q^3+1)_2(q^2-1)_2}{2^{10} } > 2 \mbox{ if } q \equiv 7 \pmod 8,\\
\frac{(q^4-1)_2(q^3+1)_2(q^2-1)_2}{2^{9} }=1 \mbox{ if } q \equiv 3 \pmod 8. \end{cases}$$

Thus, the statement $(19)$ of Theorem holds.

Finally, assume that $H$ is an almost simple maximal subgroup of $G$. Then in view of \cite[Table~8.11]{BrHoDo} $H \in \{PSL_2(7), A_7, PSL_3(4), PSU_4(2)\}$. So, $|H|_2 \le 2^6$ in view of \cite{Atlas}, and it is easy to see that $$|G|_2=\frac{(q^4-1)_2(q^3+1)_2(q^2-1)_2}{(q+1,4)}\ge 2^7.$$.

Thus, the index $|G:H|$ is even.

\medskip

Assume that $G=PSp_n(q)$, where $n \ge 4$, $n$ is even, and $q$ is odd.

 Let $H$ be the stabilizer of a non-degenerate subspace $U$ of dimension $m$ of $V$. Note that $m$ is even, and in view of (\ref{DirectSum}) we can assume that $m < n/2$ (if $m=n/2$, then $H$ is contained in a subgroup which will be considered in the next paragraph). In view of \cite[Theorem~4]{Maslova1} the index $|G:H|$ is odd if and only if $n \succeq m$. Moreover, if $n \succeq m$, then $m \not=n/2$. If  $m \not=n/2$, then in view of Tables~8.12, 8.28, 8.48, 8.64 and 8.80 of \cite{BrHoDo} and Table 3.5.C of \cite{KlLi} $H$ is maximal in $G$. Thus, the statement $(3)$ of Theorem holds for symplectic groups.

If $H$ is the stabilizer of an orthogonal decomposition $V = \bigoplus\limits_{i=1}^{n/m} V_i$ with all $V_i$ isometric  of (even) dimension $m$, then in view of \cite[Theorem~9]{Maslova1} the index $|G:H|$ is odd if and only if $m=2^w\ge 2$. In view of Tables~8.12, 8.28, 8.48, 8.64 and 8.80 of \cite{BrHoDo} and Table 3.5.C of \cite{KlLi} $H$ is maximal in $G$. Thus, the statement $(8)$ of Theorem holds.

Consider the case $n=4$. In view of the Aschbacher theorem \cite{Asch} and \cite[Tables~8.12,~8.13]{BrHoDo} any maximal subgroup of $G$ is either contained in $\bigcup_{i \in \{1,2,3,5,6\}}C_i(G)$ or is almost simple. Let $H$ be a subgroup of $G$.

Assume that $H \in C_1(G)$ and $H$ is maximal in $G$. Then  in view of \cite[Table~8.12]{BrHoDo} $H$ is the stabilizer of a totally singular subspace of $V$ of dimension $m$, where $1 \le m \le 2$. In view of \cite[Table~8.12]{BrHoDo} the preimage of $H$ in $Sp_4(q)$ is isomorphic to $E_q^{1+2}:(\mathbb{Z}_{q-1}\times Sp_2(q))$ for $m=1$ or to $E_q^3:GL_2(q)$ for $m=2$, and $$|G:H|_2=\frac{(q^4-1)_2(q^2-1)_2}{(q^2-1)_2(q-1)_2 }\ge 2.$$

Assume that  $H \in C_2(G)$. Then in view of \cite[Table~8.12]{BrHoDo} either $H$ is the stabilizer of an orthogonal decomposition $V = \bigoplus\limits_{i=1}^{2} V_i$ with $V_1$ isometric to $V_2$, and the statement $(7)$ of Theorem holds, or $H$ is the stabilizer of a decomposition $V=V_1\oplus V_2$, where $V_1$ and $V_2$ are totally singular subspaces of $V$ of dimension $2$. In the last case in view of \cite[Table~8.12]{BrHoDo} the preimage of $H$ in $Sp_4(q)$ is isomorphic to $GL_2(q).\mathbb{Z}_2$ and  $$|G:H|_2=\frac{(q^4-1)_2(q^2-1)_2}{2(q^2-1)_2(q-1)_2} \ge 2.$$

Assume that $H \in C_3(G)$. Then in view of \cite[Table~8.12]{BrHoDo} the preimage of $H$ in $Sp_4(q)$ is isomorphic either to  $Sp_2(q^2):\mathbb{Z}_2$ or to $GU_2(q).\mathbb{Z}_2$.

 In the first case

 $$|G:H|_2=\frac{(q^4-1)_2(q^2-1)_2}{2(q^4-1)_2}\ge 2.$$

In the last case

$$|G:H|_2=\frac{(q^4-1)_2(q^2-1)_2}{2(q^2-1)_2(q-1)_2} \ge 2.$$

Assume that $H \in C_5(G)$. Then $H = C_G(\sigma)$ for a field automorphism $\sigma$ of prime order $r$ of the group $G$.
If $r$ is odd, then the statement $(1)$ of Theorem holds. If $r=2$, then in view of \cite[Table~8.12]{BrHoDo} the preimage of $H$ in $Sp_4(q)$ is isomorphic to $Sp_4(q_0).\mathbb{Z}_2$ for $q=q_0^2$ and
$$|G:H|_2=\frac{(q^4-1)_2(q^2-1)_2}{2(q_0^4-1)_2(q_0^2-1)_2}= 2.$$

In view of \cite[Table~8.12]{BrHoDo} $G$ contains a maximal subgroup $H \in C_6(G)$ if and only if $q$ is prime. In view of \cite[Proposition~4.6.9]{KlLi}

$$H \cong \begin{cases} 2^4.S_5 \mbox{ if } q \equiv \pm 1 \pmod 8,\\ 2^4.A_5 \mbox{ if } q \equiv \pm 3 \pmod 8, \end{cases}$$ and

$$|G:H|_2=\begin{cases}\frac{(q^4-1)_2(q^2-1)_2}{2^{8} }\ge 2 \mbox{ if } q \equiv \pm 1 \pmod 8,\\
\frac{(q^4-1)_2(q^2-1)_2}{2^{7} }=1 \mbox{ if } q \equiv \pm 3 \pmod 8. \end{cases}$$

Thus, the statement $(21)$ of Theorem holds.

Finally, assume that $H$ is an almost simple maximal subgroup of $G$. Then in view of \cite[Table~8.13]{BrHoDo} $H$ is isomorphic either to $PSL_2(q)$, or to one of the following groups: $A_6$, $S_6$, or $A_7$. It is easy to see that the index $|G:H|$ is even.

\medskip

Assume that $G=\Omega_n(q)$, where $n \ge 7$, $n$ is odd, and $q$ is odd.

 Let $H$ be the stabilizer of a non-degenerate subspace $U$ of dimension $m$ of $V$. Note that in view of (\ref{DirectSum}) we can assume that $m$ is even. In view of \cite[Theorem~5]{Maslova1} the index $|G:H|$ is odd if and only if $n \succeq m$ and $D(U)=1$. In view of Tables~8.39, 8.58, and 8.74 of \cite{BrHoDo} and Table 3.5.D of \cite{KlLi} $H$ is maximal in $G$ if and only if $(q,m,sgn(U))\not = (3,2,+)$. However $D_2^+(3)=-1$. Thus, the statement $(4)$ of Theorem holds.

 Let $H$ be the stabilizer of an orthogonal decomposition $V = \bigoplus\limits_{i=1}^{n/m} V_i$ with all $V_i$ isometric of dimension $m$. In view of Tables~8.39, 8.58, and 8.74 of \cite{BrHoDo} and Table 3.5.D of \cite{KlLi} $H$ is maximal in $G$ if and only if either $m>1$ and $(m,q)\not = (3,3)$, or $m=1$ and $q$ is prime. In view of \cite[Theorem~10]{Maslova1} the index $|G:H|$ is odd if and only if $m=1$ and $q \equiv \pm 3 \pmod 8$. Thus, the statement $(9)$ of Theorem holds.

Consider additional abilities for the case $n=7$.

If $H \cong P\Omega_7(2)$, then in view of  \cite[Tables~8.39,~8.40]{BrHoDo} $H$ is maximal in $G$ if and only if $q$ is prime. Note that $|P\Omega_7(2)|_2=2^9$, $$|G:H|=\frac{(q^6-1)_2(q^4-1)_2(q^2-1)_2}{2^{10}},$$  and it is easy to see that the index $|G:H|$ is odd if and only if $q \equiv \pm 3 \pmod 8$. Thus, the statement $(22)$ of Theorem holds.

\medskip

Assume that $G=\Omega_n^\varepsilon(q)$, where $\varepsilon \in \{+,-\}$, $n \ge 8$, $n$ is even, and $q$ is odd.

Let $H$ be the stabilizer of a non-degenerate subspace $U$ of dimension $m$ of $V$. Note that if $m=n/2$ is odd and $U$ and $U^\bot$ are non-isometric, then in view of \cite[Theorem~12]{Maslova1} the index $|G:H|$ is even.  If $m=n/2$ and $U$ and $U^\bot$ are isometric, then $H$ is contained in a subgroup which will be considered below. Thus, if $m=n/2$, then we can assume that $m$ is even and $\varepsilon=-$. With this assumption, in view of (\ref{DirectSum}) and (\ref{MultRul2}), Tables~8.50, 8.52, 8.66, 8.68, 8.82, and 8.84 of \cite{BrHoDo} and Tables 3.5.E and 3.5.F of \cite{KlLi} $H$ is maximal in $G$ if and only if $(q,m, sgn(U)) \not = (3, 2,+), (3, n-2, \varepsilon)$. Note that $D_2^+(3)=-1$ and in view of (\ref{MultRul1}) $D_{n-2}^\varepsilon(3)=-D_n^\varepsilon(3)$.

Assume that $D(V)=1$. In view of \cite[Theorem~6]{Maslova1} the index $|G:H|$ is odd if and only if $m$ is even, $D(U)=1$, and $n \succeq m$. Note that if $n \succeq m$, then $m \not= n/2$. Thus, in the case $D(V)$=1 the statement $(5)$ of Theorem holds.

Assume that $D(V)=-1$. In view of (\ref{DirectSum}) and (\ref{MultRul1}) we can assume that $D(U)=-1$ if $m$ is even.
In view of \cite[Theorem~6]{Maslova1} the index $|G:H|$ is odd if and only if either $m$ is odd and $n-2 \succeq m-1$ (note that in this case $m \not = n/2$), or $m$ is even and $n-2 \succeq m-2$. Thus, in the case $D(V)=-1$ the statement $(5)$ of Theorem holds.

 Let $H$ be the stabilizer of an orthogonal decomposition $V = \bigoplus\limits_{i=1}^{n/m} V_i$ with all $V_i$ isometric of dimension $m$. Note that in view of \cite[Proposition~2.5.11]{KlLi} if $m$ is odd, then $D(V)=1$ (we have missed this condition in \cite[Theorem~1]{Maslova1}), and in view of (\ref{MultRul1}) and (\ref{MultRul2}) if $m$ is even, then $D(V)=(D(V_i))^{n/m}$  and $\varepsilon=(sgn(V_i))^{n/m}$. Moreover, in view of \cite[Proposition~2.5.12]{KlLi} if the mentioned conditions hold, then the corresponding decompositions exist.

 Assume that $D(V)=-1$. We conclude that $m$ is even, $n/m$ is odd, and $D(V_i)=-1$. Thus, in view of \cite[Theorem~11]{Maslova1} the index $|G:H|$ is even.

 Assume that $D(V)=1$. In view of \cite[Theorem~11]{Maslova1} the index $|G:H|$ is odd if and only if either $m=2^w\ge 2$ and $D(V_i)=1$, or $m=1$ and $q \equiv \pm 3 \pmod 8$. In view of Tables~8.50, 8.52, 8.66, 8.68, 8.82, and 8.84 of \cite{BrHoDo} and Tables 3.5.E and 3.5.F $H$ is maximal in $G$ if and only if one of the possibilities is realized: $m>1$, $(m,q) \not = (3,3)$, and $(m, q, sgn(V_i)) \not = (2, 3, \pm), (2, 5,+)${\rm;} $(n,\varepsilon) \not =(8,+)$, $m=1$, and $q$ is prime{\rm;} $(n,q)=(8,+)$, $m=1$, $q$ is prime, and $q \equiv \pm 1 \pmod 8$. Thus, the statement $(10)$ of Theorem holds.

Consider additional abilities for the case $n=8$ and $\varepsilon=+$.

If $H \cong P\Omega_8^+(2)$, then in view of  \cite[Table~8.50]{BrHoDo} $H$ is maximal in $G$ if and only if $q$ is prime. Note that $|P\Omega_8^+(2)|_2=2^{12}$, $$|G:H|=\frac{(q^4-1)(q^6-1)_2(q^4-1)_2(q^2-1)_2}{2^{14}},$$  and it is easy to see that the index $|G:H|$ is odd if and only if $q \equiv \pm 3 \pmod 8$. Thus, the statement $(23)$ of Theorem holds.  \qed\medskip

\section{Acknowledgment}

The author is very grateful to the anonymous Referee for helpful suggestions and remarks.


\end{document}